\documentclass[14pt]{article}
\usepackage{ amssymb }
\usepackage{lipsum}

\usepackage{amsfonts}
\usepackage{amsthm}

\usepackage{epsfig}
\usepackage{graphicx}  
  
\usepackage{pdfsync}

\usepackage{caption}
\usepackage{subcaption}
 
\usepackage{amsmath,amsfonts,amssymb,amsthm,epsfig,epstopdf,titling,url,array}

\usepackage{hyperref}

\theoremstyle{plain}

\newtheorem*{cor}{Corollary}

\theoremstyle{definition}

\theoremstyle{remark}

\theoremstyle{remark}

\usepackage{theoremref}

\usepackage{thmtools}

\usepackage{graphicx}
\usepackage{xcolor}
\definecolor{darkblue}{rgb}{0,0,0.5} 
\usepackage{transparent}

\DeclareMathOperator{\diam}{diam}

\newtheorem{theorem}{Theorem}[section]

\begin{document}

\title{\textbf{Logarithmically Slim Polygons}}
\author{Panagiotis Tselekidis}

\maketitle

\begin{abstract}
We study the relationship between the large-scale geometry of a geodesic
metric space and the slimness of its geodesic polygons. We show that if
there is a uniform constant $\delta$ such that every geodesic polygon in
an unbounded geodesic metric space is $\delta$-slim, then the space has
asymptotic Assouad--Nagata dimension one.

We further consider graphs of bounded valency and allow the slimness
constant of an $n$-gon to grow logarithmically with $n$. More precisely,
under an explicit logarithmic bound on the slimness of geodesic
$n$-gons, we prove that the graph has asymptotic Assouad--Nagata
dimension one. Thus, sufficiently strong logarithmic control on the
slimness of geodesic polygons forces one-dimensional large-scale
geometry.

We also obtain two criteria for graphs, outside the hyperbolic
setting, to have asymptotic Assouad--Nagata dimension one. In particular,
these criteria show that one-dimensional large-scale geometry can be
detected from quantitative restrictions on the slimness of geodesic
polygons, without requiring the space to be hyperbolic.


\end{abstract}

\section{Introduction}

The geometry of geodesic polygons provides a useful way of measuring
how far a metric space is from being tree-like. In a tree, every
geodesic polygon is degenerate in the sense that each side is contained
in the union of the remaining sides. More generally, hyperbolic spaces
are characterized by the slimness of geodesic triangles. This raises a
natural question: \begin{center}
\emph{what can be deduced if one imposes quantitative
slimness conditions on geodesic polygons with an arbitrary number of
sides?}
\end{center}

In this paper we investigate this question from the point of view of
asymptotic Assouad--Nagata dimension.

Recall that the asymptotic
Assouad--Nagata dimension, denoted by $\operatorname{asdim}_{AN}$,
is a large-scale dimension which strengthens the usual asymptotic
dimension by requiring linear control of the diameters of the members
of the covers. Thus, while asymptotic dimension is a coarse invariant
measuring the large-scale dimension of a space, asymptotic
Assouad--Nagata dimension retains quantitative information about the
scale at which this dimension is witnessed.\\

Our first result concerns spaces in which the slimness of geodesic
polygons is uniformly bounded, independently of the number of sides.

\begin{theorem}
Let $X$ be an unbounded geodesic metric space. Suppose that there exists
$\delta \geq 0$ such that every geodesic $n$-gon in $X$ is
$\delta$-slim. Then
\[
\operatorname{asdim}_{AN} X = 1.
\]
\end{theorem}

The hypothesis is substantially stronger than ordinary Gromov
hyperbolicity. 
The conclusion reflects this stronger form of thinness: the
large-scale geometry is one-dimensional. In particular,
$\operatorname{asdim}X=1$, and, by the characterization of
hyperbolic spaces of asymptotic dimension one, $X$ is quasi-isometric to a tree.

We observe that in every $\delta$-hyperbolic geodesic metric space, every geodesic $n$-gon is $(n-2)\delta$-slim. This raises the following natural question:

\begin{flushleft}
\textbf{Question 1:}\label{question1}
Is there an increasing function $f\colon (0,\infty) \times \{\delta\}\to[0,\infty)$ such that $\lim_n f(n,\delta)= \infty$, and  if every geodesic $n$-gon in a graph $G$ is $f(n,\delta)$-slim, then $
\operatorname{asdim}_{AN}G=1$?
\end{flushleft}

The main point of the paper is to answer question \ref{question1}. We show that, for graphs of
bounded valency, one may allow the slimness constant to grow
logarithmically with the number of sides.

Let $G$ be a graph of valency at most $k$, equipped with its path
metric. We say that a geodesic $n$-gon is $\frac{1}{8}\log_k(n\delta)\text{-slim}$ if each of its sides is contained in the $\frac{1}{8}\log_k(n\delta)$-neighbourhood of the union of the remaining sides. Our second main result is the following.

\begin{theorem}
Let $G$ be an unbounded connected graph of valency at most $k$.
Suppose that there exists $\delta>0$ such that every geodesic
$n$-gon in $G$ is
\[
\frac{1}{8}\log_k(n\delta)\text{-slim}.
\]
Then
\[
\operatorname{asdim}_{AN}G=1.
\]
\end{theorem}


Our result
shows that, below the quantitative threshold appearing in the
hypothesis, this restriction is sufficient to force linear control and
hence asymptotic Assouad--Nagata dimension one.\\

As a corollary we obtain the following characterisation of quasi-trees with bounded valency:

\begin{cor}\label{Cor1}
Let $X$ be an unbounded connected graph of valency at most $k$.\\
Then, there exists $\delta>0$ such that every geodesic 
$n$-gon in $G$ is $\frac{1}{8}\log_k(n\delta)\text{-slim}$ if and only if,
$X$ is quasi-isometric to a tree.

\end{cor}



In Euclidean spaces, geodesic polygons satisfy explicit
quantitative bounds on their slimness. More precisely, if $P$ is a
geodesic $n$-gon in a Euclidean space and
$l_0\leq l_1\leq\cdots\leq l_n$ are the side lengths of $P$, then $P$ is  

\begin{center}
$\frac{n-2}{2} l_{n-2}$-slim
\end{center}

This naturally leads to the question of whether analogous quantitative
restrictions on polygonal slimness can force one-dimensional
large-scale geometry in more general geodesic spaces.

We show the following:

\begin{theorem}
Let $X$ be a geodesic graph with valency at most $k$. Suppose that every geodesic $n$-gon, with side lengths $l_0\leq l_1\leq\cdots\leq l_{n}, $ is 

$$\frac{1}{8}\log_k\left((n-2)\frac{l_{n-2}}{17}\right)
\text{-slim}.$$ Then

$$
\operatorname{asdim}_{AN}X=1.
$$

\end{theorem}

For vertex transitive graphs we have a better bound.

\begin{theorem}
Let $X$ be a vertex-transitive geodesic graph of valency at most $k$. Suppose that every geodesic $n$-gon with side lengths
 $l_0 \leq l_1 \leq \dots \leq l_n$

is

$$
\frac{1}{8}\log_k\left(\frac{(n-2)l_{n-2}}{2\varepsilon}\right)\text{-slim}
$$

for some $\varepsilon>1$. Then

$$
\operatorname{asdim}_{AN}X=1.
$$

\end{theorem}

These two criteria show that quantitative restrictions on the slimness
of geodesic polygons can force one-dimensional large-scale geometry
even outside the hyperbolic setting.\\

The paper is organized as follows. In Section~\ref{prelims}
we recall the definitions and basic properties of asymptotic
Assouad--Nagata dimension and slim geodesic polygons. In
Section~\ref{section3} we prove the result for uniformly slim
polygons. In Section~\ref{section4} we establish the logarithmic
polygon slimness theorem for bounded-valency graphs.
In Section~\ref{section5}, we prove the last two criteria for more general graphs, without assuming hyperbolicity.
\section{Preliminaries}\label{prelims}

If $(X,d)$ is a metric space and $A\subseteq X$, we write
$N_r(A)=\{x\in X:d(x,A)\le r\}$, for the closed $r$-neighbourhood of $A$.\\
Let \(X\) be a metric space and let \(R>0\). An \emph{\(R\)-path} in \(X\) is a finite sequence of points $ x_0,x_1,\ldots,x_m\in X $ such that $ d(x_{i-1},x_i)\le R$ for every $i=1,\ldots,m. $ We say that a subset $Z$ of $X$ is $R-$connected id every two points in $Z$ are connected by an $r$-path. Let $C \subseteq X$, we say that $C$ is an \emph{$R-$component} of $X$ if it is $R$-connected, and for every other $R$-connected subset $Z$, of $X$ such that $C \subseteq Z$, we have that $C=Z$.

The \emph{length} of an \(R\)-path is defined by

$$
\ell(x_0,\ldots,x_m)
=\sum_{i=1}^{m}d(x_{i-1},x_i).
$$

The \emph{\(R\)-length} of the path is the number of \(R\)-edges, namely

$$
\ell_R(x_0,\ldots,x_m)=m.
$$

We say that the \(R\)-path joins \(x_0\) to \(x_m\).

\subsection{Geodesic metric spaces and graphs}

A metric space $(X,d)$ is called \emph{geodesic} if for every $x,y\in X$ there exists an isometric embedding $ \gamma\colon [0,d(x,y)]\longrightarrow X$, such that $\gamma(0)=x$ and $\gamma(d(x,y))=y$.
The image of such an embedding is called a \emph{geodesic} from $x$ to $y$ and is denoted by $[x,y]$. A connected graph will always be regarded as a metric space by giving every edge length one and considering the induced path metric. In particular, graphs are geodesic metric spaces.

For a graph $G$, we denote its valency by

$$
\operatorname{val}(G)=\sup_{v\in V(G)}\deg(v).
$$

We say that $G$ has \emph{bounded valency} if $\operatorname{val}(G)<\infty$.

\subsection{Geodesic polygons and slimness}

Let $X$ be a geodesic metric space. A \emph{geodesic $n$-gon} is a cyclic sequence

$$
P=[x_1,x_2]\cup [x_2,x_3]\cup\cdots\cup[x_n,x_1],
$$

where each $[x_i,x_{i+1}]$ is a geodesic.

For $\delta\ge0$, a geodesic $n$-gon $P$ is called \emph{$\delta$-slim} if every side of $P$ is contained in the $\delta$-neighbourhood of the union of the remaining sides. Thus, for every $i$, $[x_i,x_{i+1}]
\subseteq
N_\delta\left(
\bigcup_{j\ne i}[x_j,x_{j+1}]
\right).
$

A geodesic metric space is called \emph{$\delta$-hyperbolic} if every geodesic triangle is $\delta$-slim. It is called \emph{hyperbolic} if it is $\delta$-hyperbolic for some $\delta\ge0$.

Notice that hyperbolicity only imposes a uniform bound on geodesic triangles. The corresponding bound for general geodesic polygons need not be uniform in the number of sides. Indeed, a standard consequence of hyperbolicity is that an $n$-gon is $\varepsilon_n$-slim, where $\varepsilon_n$ depends on $n$ and $\delta$, where $
\varepsilon_n\le (n-2)\delta.
$



\subsection{Asymptotic Assouad--Nagata dimension}

We recall Gromov's asymptotic dimension. The \emph{asymptotic dimension} $\operatorname{asdim} X$ of a metric space $X$ is defined as follows. We say that $\operatorname{asdim} X\le n$
if and only if, for every $R>0$, there exists a uniformly bounded cover $\mathcal U$ of $X$ whose $R$-multiplicity is at most $n+1$; that is, every ball of radius $R$ in $X$ intersects at most $n+1$ elements of $\mathcal U$.

The \emph{asymptotic Assouad--Nagata dimension} $\operatorname{asdim}_{AN} X$ of a metric space $X$ is defined similarly, with the additional requirement of linear control on the diameters of the elements of the cover. More precisely,

$$
\operatorname{asdim}_{AN} X\le n
$$

if and only if there exist constants $C, R_0 >0$ such that, for every $R>R_0$, there exists a cover $\mathcal U$ of $X$ with $R$-multiplicity at most $n+1$ and

$$
\operatorname{diam}(U)\le CR
$$

for every $U\in\mathcal U$. It is easy to see that $ \operatorname{asdim} X \leq \operatorname{asdim}_{AN} X$. The terminology and basic theory of asymptotic Assouad--Nagata dimension can be found in \cite{DS07}.

\section{Uniformly slim polygons}\label{section3}
\begin{theorem}
Let $X$ be an unbounded geodesic metric space. Suppose that there exists
$\delta \geq 0$ such that every geodesic $n$-gon in $X$ is
$\delta$-slim. Then
\[
\operatorname{asdim}_{AN} X = 1.
\]
\end{theorem}
\begin{proof}
Let $R>10\delta$, we will construct a suitable uniformly bounded cover for $X$, say $\mathcal{U}$.

Let $e\in X$ (base point). We decompose $X$ using annuli $A_i=B_{3R(i+1)}(e)\setminus B_{3Ri}(e)$, where $i \in \mathbb{N}$. We will construct suitable covers, say $\mathcal{U}_i$, for all annuli $A_i$, for $i>\geq 2$. 
For the annulus $A_i$, we consider its R-components.\\ \textbf{Claim:} If $C$ is an R-component of $A_i$ ($i\geq2$), then $\operatorname{diam}C<8 R.$
\begin{proof}
Let $x,y \in C$, and let $x=x_0, \dots, x_n=y$ be an R-path from $x$ to $y$. Let $r_x$, $r_y$ be geodesics from $e$ to $x$, and from $e$ to $y$, respectively. Observe that the geodesics $r_x, [x_0,x_1], \dots, [x_{n-1},x_n], r_y $ form a geodesic $(n+2)-$gon. From the hypothesis this polygon is $\delta-$slim, so $x'=r_x(3Ri-(\frac{R}{2}+2\delta))$ is $\delta$-close to one of the other geodesics.\\
However, from the way we choose $x'$, we can see that $\operatorname{dist}(x',[x_j,x_{j+1}])>\delta$. Thus, $\operatorname{dist}(x',r_y)\leq\delta$. In particular, we can deduce that $d(x',y')\leq 2\delta$, where $y'=r_y(3Ri-(\frac{R}{2}+2\delta))$. Then $d(x,y)\leq 3R+(\frac{R}{2}+2\delta)+2\delta+ 3R+(\frac{R}{2}+2\delta)< 8R.$
\end{proof}  
We set $\mathcal{U}_i$, to be the collection of the R-components of the annulus $A_i$, $i\geq2.$ Then the set $\mathcal{U}= B_{6R}(e) \cup \bigcup_{i\geq 2}\mathcal{U}_i$, is a cover of $X$, with R-multipicity=2 ($X$ is unbounded and geodesic), and every element in $\mathcal{U}$ has diameter at most $12R$. It follows that $\operatorname{asdim}_{AN} X = 1.$
\end{proof}

Since $X$ is unbounded and $\operatorname{asdim} X \leq \operatorname{asdim}_{AN} X=1 $, we can deduce as a corollary that: 

\begin{cor}\label{Cor1}
Let $X$ be an unbounded geodesic metric space. Suppose that there exists
$\delta \geq 0$ such that every geodesic $n$-gon in $X$ is
$\delta$-slim. Then $X$ is quasi-isometric to a tree.
\end{cor}

\section{Logarithmically slim polygons}\label{section4}
In this section we will give an answer to Question \ref{question1}.

\begin{theorem}\label{thm4.1}
Let $X$ be an unbounded connected graph with $\operatorname{val}(X)=k$.
Suppose that there exists $\delta>0$ such that every geodesic
$n$-gon in $G$ is $\frac{1}{8}\log_k(n\delta)\text{-slim}$.
Then
\[
\operatorname{asdim}_{AN}G=1.
\]
\end{theorem}
\begin{proof}
Let $R>100(\delta+1)$, we will construct a suitable uniformly bounded cover for $X$, say $\mathcal{U}$.

Let $e\in X$ (base point). We decompose $X$ using annuli $A_i=B_{3R(i+1)}(e)\setminus B_{3Ri}(e)$, where $i \in \mathbb{N}$. We will construct suitable covers, say $\mathcal{U}_i$, for all annuli $A_i$, for $i\geq 17$. 
For the annulus $A_i$, we consider its R-components.\\ \textbf{Claim:} If $C$ is an R-component of $A_i$ ($i\geq 17$), then $\operatorname{diam}C \leq 45R.$
\begin{proof}
Suppose for a contradiction that $\diam C >mR$, where $m=45$. Let $x,y \in C$, and let $x=x_0, \dots, x_n=y$ be an R-path from $x$ to $y$ of minimal R-length, we denote this path by $\gamma$. There exists $ 0<\epsilon<1$ such that 
\begin{equation}\label{eq1}
d(x,y)=mR+ \epsilon R.
\end{equation}  Let $r_x$, $r_y$ be geodesics from $e$ to $x$, and from $e$ to $y$, respectively. Observe that the geodesics $r_x, [x_0,x_1], \dots, [x_{n-1},x_n], r_y $ form a geodesic $(n+2)-$gon. From the hypothesis, this polygon is $\frac{1}{8}\log_k((n+2)\delta)\text{-slim}$, so $x'=r_x(3Ri-(\frac{3R}{5}+\frac{1}{8}\log_k((n+2)\delta))$ is $\frac{1}{8}\log_k((n+2)\delta)$-close to one of the other geodesics. 

As we see in the case 1 and 2 bellow, we can choose suitable $x , y$ such that $\gamma \subseteq B_{d(x,y)}(x)$.
Then it follows that $\log_k ((n+2)\delta) \leq 47R$. Thus for $i \geq 17$, we have that $3Ri>\frac{3R}{5}+\frac{1}{8}\log_k((n+2)\delta)$.

From the way we choose $x'$, we can see that $\operatorname{dist}(x',[x_j,x_{j+1}])>\frac{1}{8}\log_k((n+2)\delta)$. Thus, $\operatorname{dist}(x',r_y) \leq \frac{1}{8}\log_k((n+2)\delta)$.
In particular, we can deduce that $d(x',y')\leq 2\frac{1}{8}\log_k((n+2)\delta)$, where $y'=r_y(3Ri-(\frac{3R}{5}+\frac{1}{8}\log_k((n+2)\delta)))$. 

Let $x''=r_x(3Ri)$ and $y''=r_y(3Ri)$. \emph{Since both $x'', y'' \in C$, to prove the claim it is enough to show that $d(x'', y'')<38R$.} So we assume that $x=x''$ and $y=y''$. So we have $d(x,y)\leq \frac{6R}{5} + \frac{1}{2}\log_k((n+2)\delta))$, we deduce that

\begin{equation}\label{eq2}
(m-\frac{6}{5})R\leq  \frac{1}{2}\log_k(n(2\delta))
\end{equation}

We set $\overline{\gamma}$ to be the path formed by a union of the geodesics  $ [x_0,x_1], \dots, [x_{n-1},x_n] $.\\ 
\textbf{Observation 1}\label{ob1}
Since $\gamma$ has minimal R-length we observe that the set $\{ x_{2k} \mid k \in \mathbb{N}\}$ is $\frac{R}{2}$-disjoint.\\

We distinguish two cases:
\begin{flushleft}
\textbf{Case 1:} $\gamma \subseteq B_{d(x,y)}(x)$.\\
Then $\overline{\gamma} \subseteq B_{d(x,y)+\frac{R}{2}}(x)$. Thus by observation 1 there are at least $\frac{n}{2}-1$ many disjoint $\frac{R}{2}-$balls in $B_{d(x,y)+\frac{R}{2}}(x).$ Thus, there are at least $(\frac{n}{2}-1)\frac{R}{2}$ many elements in $B_{d(x,y)+\frac{R}{2}}(x).$ So,  

$$(\frac{n}{2}-1)\frac{R}{2} \leq \mid B_{d(x,y)+\frac{R}{2}}(x)\mid \leq k^{d(x,y)+\frac{R}{2}}+1.$$

If $n<4$, then $d(x,y)<4R$. If both $n\geq4$, then   

$$\frac{nR}{8}  \leq 2k^{(m+2)R}.$$

Then, 

\begin{equation}\label{eq3}
log_k(\frac{nR}{16}) \leq (m+2)R
\end{equation}


Since $m> 5$. We have that $m+2\leq \frac{7m}{5}$, and $\frac{5m}{7} \leq m - \frac{6}{5}$. Thus inequalities \ref{eq2} and \ref{eq3} will transform to the following inequalities: 

\begin{equation}\label{eq4}
\frac{5m}{7}R\leq   \frac{1}{2}\log_k(n(2\delta))
\end{equation}

\begin{equation}\label{eq5}
log_k(\frac{nR}{16}) \leq \frac{7m}{5}R
\end{equation}

Combining them we obtain the following: 

\begin{equation}\label{eq6}
log_k(\frac{nR}{16}) \leq  \frac{7^2}{5^2}\frac{1}{2}\log_k(n(2\delta)).
\end{equation}
Since $R>32\delta$, and $\frac{7^2}{5^2}\frac{1}{2}<1$, we have found a contradiction. 

\end{flushleft}

\begin{flushleft}
\textbf{Case 2:} $\gamma \nsubseteq B_{d(x,y)}(x)$.\\
Then there exists $n_y \in \{ 1, \dots, n \} $ such that $x_j \in B_{d(x,y)}(x)$, for all $j \leq n_y$, while $x_{n_y +1} \notin B_{d(x,y)}(x)$. Then $d(x,y''') \in [d(x,y)-R, d(x,y)]$, where $y'''=x_{n_y }$.\\
Observe that the R-path $x=x_0, \dots,x_{n_y }=y'''$ belongs to $B_{d(x,y)} \cap A_i$, so we can use the arguments (slightly modified on constants) from Case 1. So the analogue of inequality \ref{eq2} is the following
 
 \begin{equation}\label{eq7}
(m-\frac{26}{5})R\leq  \frac{1}{2}\log_k(n_y(2\delta))
\end{equation}
 
 and the analogue of inequality \ref{eq3} is 
 
 \begin{equation}\label{eq8}
log_k(\frac{n_yR}{16}) \leq (m+2)R
\end{equation}
 
 Since $m >37$. We have $m+2\leq \frac{7m}{5}$, and $\frac{5m}{7} \leq m - \frac{26}{5}$. Thus inequalities \ref{eq7} and \ref{eq8} will transform to the following inequalities: 

\begin{equation}\label{eq4}
\frac{5m}{7}R\leq   \frac{1}{2}\log_k(n_y(2\delta))
\end{equation}

\begin{equation}\label{eq5}
log_k(\frac{n_yR}{16}) \leq \frac{7m}{5}R
\end{equation}

Combining them we obtain the following: 

\begin{equation}\label{eq6}
log_k(\frac{n_yR}{16}) \leq  \frac{7^2}{5^2}\frac{1}{2}\log_k(n_y(2\delta)).
\end{equation}
Since $R>32\delta$, and $\frac{7^2}{5^2}\frac{1}{2}<1$, we have found a contradiction.

\end{flushleft}

\end{proof}

We set $\mathcal{U}_i$, to be the collection of the R-components of the annulus $A_i$, $i\geq 17.$ Then the set $\mathcal{U}= B_{51R}(e) \cup \bigcup_{i\geq 17}\mathcal{U}_i$, is a cover of $X$, with R-multipicity=2 ($X$ is unbounded and geodesic), and every element in $\mathcal{U}$ has diameter at most $102R$. It follows that $\operatorname{asdim}_{AN} X = 1.$

\end{proof}

As a corollary we obtain the following characterisation of quasi-trees with uniformly bounded valency:

\begin{cor}\label{Cor1}
Let $X$ be an infinite connected graph with $\operatorname{val}(X)=k$.\\
Then, there exists $\delta>0$ such that every geodesic 
$n$-gon in $G$ is $\frac{1}{8}\log_k(n\delta)\text{-slim}$ if and only if,
$X$ is quasi-isometric to a tree.

\end{cor}

From the previous proof we can observe that we can get a slightly closer to 1 multiplicity constant, instead of $\frac{1}{8}$, on the slimness of polygons, however, the following questions are open:

\begin{flushleft}
\textbf{Question 2:}\label{question2}
Does the analogue of Theorem \ref{thm4.1} hold if every geodesic $n$-gon in $X$ is $log_k(n\delta)\text{-slim}$? 
\end{flushleft}

We can go a bit further, and ask:

\begin{flushleft}
\textbf{Question 3:}\label{question3}
Does the analogue of Theorem \ref{thm4.1} hold if every geodesic $n$-gon in $X$ is $(log_k(n\delta))^2\text{-slim}$? 
\end{flushleft}

\begin{flushleft}
\textbf{Question 4:}\label{question3}
Does the analogue of Theorem \ref{thm4.1} hold if every geodesic $n$-gon in $X$ is $ \sqrt{n\delta}\text{-slim}$? 
\end{flushleft}

\section{Logarithmically slim polygons on general graphs}\label{section5}
In Euclidean spaces, geodesic polygons satisfy explicit
quantitative bounds on their slimness. More precisely, if $P$ is a
geodesic $n$-gon in a Euclidean space and
$l_0\leq l_1\leq\cdots\leq l_n$ are the side lengths of $P$, then $P$ is  

\begin{center}
$(n-2)l_{n-2}$-slim
\end{center}

Indeed, it is not hard to prove it for triangles. Then using induction on the number of vertices of the n-gons, and the fact that we can separate the two largest sides of an n-gon P by a geodesic we can obtain the wanted bound.

The above bound leads to the question of whether analogous quantitative
restrictions on polygonal slimness can force one-dimensional
large-scale geometry in more general geodesic spaces. To answer this question we show the following:

\begin{theorem}\label{thm5.1}
Let $X$ be an infinite graph with $\operatorname{val}(X)=k$. Suppose that every geodesic $n$-gon, with side lengths $l_0\leq l_1\leq\cdots\leq l_{n}, $ is 

$$\frac{1}{8}\log_k\left((n-2)\frac{l_{n-2}}{17}\right)
\text{-slim}.$$ Then

$$
\operatorname{asdim}_{AN}X=1.
$$

\end{theorem}

For vertex transitive graphs we have a better bound.

\begin{theorem}\label{thm5.2}
Let $X$ be a vertex-transitive infinite graph with $\operatorname{val}(X)=k$. Suppose that every geodesic $n$-gon with side lengths
 $l_0 \leq l_1 \leq \dots \leq l_n$

is

$$
\frac{1}{8}\log_k\left(\frac{(n-2)l_{n-2}}{2\varepsilon}\right)\text{-slim}
$$

for some $\varepsilon>1$. Then

$$
\operatorname{asdim}_{AN}X=1.
$$

\end{theorem}
Since the proofs of the two theorems are similar, we give only the proof of the first, pointing out along the way any steps that require a different argument for the second.

\begin{proof}

Let $R>100(\delta+1)$, we will construct a suitable uniformly bounded cover for $X$, say $\mathcal{U}$.

Let $e\in X$ (base point). We decompose $X$ using annuli $A_i=B_{3R(i+1)}(e)\setminus B_{3Ri}(e)$, where $i \in \mathbb{N}$. We will construct suitable covers, say $\mathcal{U}_i$, for all annuli $A_i$, for $i\geq 17$. 
For the annulus $A_i$, we consider its R-components.\\ \textbf{Claim:} If $C$ is an R-component of $A_i$ ($i\geq 17$), then $\operatorname{diam}C \leq 45R.$
\begin{proof}
Suppose for a contradiction that $\diam C >mR$, where $m=45$. Let $x,y \in C$, and let $x=x_0, \dots, x_n=y$ be an R-path from $x$ to $y$ of minimal R-length, we denote this path by $\gamma$. There exists  $ 0<\epsilon<1$ such that 
\begin{equation}\label{eq1}
d(x,y)=mR+ \epsilon R.
\end{equation}  Let $r_x$, $r_y$ be geodesics from $e$ to $x$, and from $e$ to $y$, respectively. Observe that the geodesics $r_x, [x_0,x_1], \dots, [x_{n-1},x_n], r_y $ form a geodesic $(n+2)-$gon. From the hypothesis, this polygon is $\frac{1}{8}\log_k\left(n\frac{l_{n}}{17}\right)
\text{-slim}.$, so $x'=r_x(3Ri-(\frac{3R}{5}+\frac{1}{8}\log_k\left(n\frac{l_{n}}{17}\right)$ is $\frac{1}{8}\log_k\left(n\frac{l_{n}}{17}\right)$-close to one of the other geodesics.

As we will see in the case 1 and 2 bellow, we can choose suitable $x , y$ such that $\gamma \subseteq B_{d(x,y)}(x)$.
Then it follows that $\log_k (n l_{n}) \leq 48R$. Thus for $i \geq 17$, we have that $3Ri>\frac{3R}{5} + \frac{1}{8} \log_k(n \frac{l_{n}}{17}), \frac{3R}{5} + \frac{1}{8} \log_k(n \frac{l_{n}}{2 \varepsilon})$.

From the way we choose $x'$, we can see that $\operatorname{dist}(x',[x_j,x_{j+1}])>\frac{1}{8}\log_k\left(n\frac{l_{n}}{17}\right)$. Thus, $\operatorname{dist}(x',r_y) \leq \frac{1}{8}\log_k\left(n\frac{l_{n}}{17}\right)$.
In particular, we can deduce that $d(x',y')\leq 2\frac{1}{8}\log_k\left(n\frac{l_{n}}{17}\right)$, where $y'=r_y(3Ri-(\frac{3R}{5}+\frac{1}{8}\log_k\left(n\frac{l_{n}}{17}\right))$. 

We have $d(x,y)\leq 9R + \frac{1}{2}\log_k\left(n\frac{l_{n}}{17}\right)$, we deduce that

\begin{equation}\label{eq13}
(m-9)R\leq  \frac{1}{2}\log_k\left(n\frac{l_{n}}{17}\right)
\end{equation}

\emph{In the proof of the second theorem we have }
\begin{equation}\label{eq14}
(m-9)R\leq  \frac{1}{2}\log_k\left(n\frac{l_{n}}{2\varepsilon}\right)
\end{equation}

We set $\overline{\gamma}$ to be the path formed by a union of the geodesics  $ [x_0,x_1], \dots, [x_{n-1},x_n] $.\\ 
\textbf{Observation 1}\label{ob1}
Since $\gamma$ has minimal R-length we observe that the set $\{ x_{2k} \mid k \in \mathbb{N}\}$ is $\frac{R}{2}$-disjoint.\\

We distinguish two cases:
\begin{flushleft}
\textbf{Case 1:} $\gamma \subseteq B_{d(x,y)}(x)$.\\
Then $\overline{\gamma} \subseteq B_{d(x,y)+\frac{R}{2}}(x)$. Thus by observation 1 there are at least $\frac{n}{2}-1$ many disjoint $\frac{R}{2}-$balls in $B_{d(x,y)+\frac{R}{2}}(x).$ Thus, there are at least $(\frac{n}{2}-1)\frac{R}{2}$ many elements in $B_{d(x,y)+\frac{R}{2}}(x).$ So,  

\begin{equation}\label{eq15}
(\frac{n}{2}-1)\frac{R}{2} \leq \mid B_{d(x,y)+\frac{R}{2}}(x)\mid \leq k^{d(x,y)+\frac{R}{2}}+1.
\end{equation}

If $n<4$, then $d(x,y)<4R$. If $n\geq4$, then   

$$\frac{nR}{8}  \leq 2k^{(m+2)R}.$$

Then, 

\begin{equation}\label{eq16}
log_k(\frac{nR}{16}) \leq (m+2)R
\end{equation}


Since $m> 10$. We have that $m+2\leq \frac{7m}{5}$, and $\frac{5m}{7} \leq m - 9$. Thus inequalities \ref{eq13} and \ref{eq16} will transform to the following inequalities: 

\begin{equation}\label{eq17}
\frac{5m}{7}R\leq   \frac{1}{2}\log_k\left(n\frac{l_{n}}{17}\right)
\end{equation}

\begin{equation}\label{eq18}
log_k(\frac{nR}{16}) \leq \frac{7m}{5}R
\end{equation}

Combining them we obtain the following: 

\begin{equation}\label{eq19}
log_k(\frac{nR}{16}) \leq  \frac{7^2}{5^2}\frac{1}{2}\log_k\left(n\frac{l_{n}}{17}\right).
\end{equation}
Here we make they \emph{key observation} that $l_n$ is the length of one of the geodesics $[x_j,x_{j+1}]$. Indeed, the geodesics $r_x$, $r_y$ have length at least $6R$  (so $\{l(r_x)$, $l(r_y)\}= \{l_{n+1},l_{n+2}\})$), while the geodesics $[x_j,x_{j+1}]$ have length at most $R$.   
Since $R\leq l_n$, and $\frac{7^2}{5^2}\frac{1}{2}<1$, we have found a contradiction. 

\textbf{For the proof of the second theorem:}
We assume that $R>log_k(\frac{1}{\sqrt{\varepsilon}-1}) $.
By observation 1 there are at least $\frac{n}{2}-1$ many disjoint $\frac{R}{2}-$balls in $B_{d(x,y)+\frac{R}{2}}(x).$
We observe that since $X$ is vertex transitive with bounded valency there exists a bi-infinite geodesic passing from every vertex of $X$. Thus, there are at least $(\frac{n}{2}-1)R$ many elements in $B_{d(x,y)+\frac{R}{2}}(x).$ So, 

\begin{equation}\label{eq20}
(\frac{n}{2}-1)R\leq  k^{d(x,y)+\frac{R}{2}}+1.
\end{equation}

If  $n\leq \frac{2\sqrt{\varepsilon} }{ \sqrt{\varepsilon}-1}$, then $d(x,y)<  (\frac{2\sqrt{\varepsilon} }{\sqrt{\varepsilon}-1 }+1)R$. If  $n>\frac{2\sqrt{\varepsilon} }{\sqrt{\varepsilon}-1 }$, then   

$$\frac{nR}{2 \sqrt{\varepsilon}}  \leq \sqrt{\varepsilon}k^{(m+2)R}.$$

Then, 

\begin{equation}\label{eq21}
log_k(\frac{nR}{2\varepsilon}) \leq (m+2)R
\end{equation}


Since $m > 10$. We have that $m+2\leq \frac{7m}{5}$, and $\frac{5m}{7} \leq m - 9$. Thus inequalities \ref{eq13} and \ref{eq21} will transform to the following inequalities: 

\begin{equation}\label{eq17}
\frac{5m}{7}R\leq   \frac{1}{2}\log_k\left(n\frac{l_{n}}{2\varepsilon}\right)
\end{equation}

\begin{equation}\label{eq18}
log_k(\frac{nR}{2\varepsilon}) \leq \frac{7m}{5}R
\end{equation}

Combining them we obtain the following: 

\begin{equation}\label{eq19}
log_k(\frac{nR}{2\varepsilon}) \leq  \frac{7^2}{5^2}\frac{1}{2}\log_k\left(n\frac{l_{n}}{2\varepsilon}\right).
\end{equation}
We observe similarly, that $l_n$ is the length of one of the geodesics $[x_j,x_{j+1}]$.  
Since $R\leq l_n$, and $\frac{7^2}{5^2}\frac{1}{2}<1$, we have found a contradiction. So, $d(x,y)\leq \max \{45, \frac{2\sqrt{\varepsilon} }{ \sqrt{\varepsilon}-1} +1 \}R$.

\end{flushleft}

\begin{flushleft}
\textbf{Case 2:} $\gamma \nsubseteq B_{d(x,y)}(x)$.\\
Then there exists $n_y \in \{ 1, \dots, n \} $ such that $x_j \in B_{d(x,y)}(x)$, for all $j \leq n_y$, while $x_{n_y +1} \notin B_{d(x,y)}(x)$. Then $d(x,y''') \in [d(x,y)-R, d(x,y)]$, where $y'''=x_{n_y }$.\\
Observe that the R-path $x=x_0, \dots,x_{n_y }=y'''$ belongs to $B_{d(x,y)} \cap A_i$, so we can use the arguments for $x$ and $y'''$ from Case 1. In the end we will obtain that $d(x,y''')<cR$, where $c=\max \{11, \max \{45, \frac{2\sqrt{\varepsilon} }{ \sqrt{\varepsilon}-1} +1 \}\}$. So, $d(x,y)<\leq cR$.

\end{flushleft}

\end{proof}
We set $\mathcal{U}_i$, to be the collection of the R-components of the annulus $A_i$, $i\geq17.$ Then the set $\mathcal{U}= B_{51R}(e) \cup \bigcup_{i\geq 17}\mathcal{U}_i$, is a cover of $X$, with R-multipicity=2 ($X$ is unbounded and geodesic), and every element in $\mathcal{U}$ has diameter at most $102R$ for theorem \ref{thm5.1}, and at most $\max \{102, \frac{2\sqrt{\varepsilon} }{ \sqrt{\varepsilon}-1} +1 \}R$ for theorem \ref{thm5.2}. It follows that $\operatorname{asdim}_{AN} X = 1.$

\end{proof}

\begin{flushleft}
Department of Mathematics University of the Aegean Karlovassi, 83200 Samos, Greece.\\
\textit{E-mail address}: ptselekidis@aegean.gr
\end{flushleft}

\end{document}